%% file: non-cox_e.tex
\documentclass{article}
\usepackage{latexsym,amsthm,amsfonts,epsfig,psfrag }
\begin{document}

\newtheorem{cor}{Corollary}
\newtheorem{lemma}{Lemma}
\newtheorem{theorem}{Theorem}
\newtheorem{prop}{Proposition}
\renewcommand{\proofname}{Proof}

\def \Z{{\mathbb Z}}
\def \S{{\mathbb S}}
\def \H{{\mathbb H}}
\def \E{{\mathbb E}}


\begin{center}
{\Large \bf
Euclidean simplices generating \\
discrete reflection groups}

\vspace{7pt}
{\large A.~Felikson, P.~Tumarkin }
\end{center}

\section*{Introduction}

\noindent

Let $P$ be a convex polytope in  the spherical space $\S^n$,
in the Euclidean space $\E^n$, or in the hyperbolic space $\H^n$.
 Consider the group $G_P$ generated by
 reflections in the facets of $P$.
We call $G_P$  a {\it reflection group generated by $P$}.
The problem we consider in this paper
is to list polytopes generating discrete reflection groups.

The answer is known only for
some  combinatorial types of polytopes.
Already in 1873, Schwarz \cite{Shwarz}  listed spherical triangles
generating discrete groups.
In  1998, E.~Klimenko and M.~Sakuma \cite{Pink}
solved the problem
for hyperbolic triangles.
In \cite{polygons}, \cite{pyrprism}, \cite{simplex},
\cite{lambert} 
the problem was solved for hyperbolic quadrilaterals,
compact hyperbolic pyramids and triangular prisms, hyperbolic simplices,
and Lambert cubes in $\S^3$, $\E^3$, $\H^3$.
In~\cite{sph} the problem was solved for spherical simplices. 

\vspace{5pt}
\noindent

In this paper, we use the method of~\cite{sph} to classify Euclidean
simplices generating discrete reflection groups.

\vspace{5pt}
\noindent

The authors are grateful to the Max-Planck Institute for Mathematics
in Bonn for hospitality.

\vspace{20pt}

\section{Preliminaries}

A convex polytope in $\S^n$, $\E^n$ or $\H^n$ is called a 
{\it Coxeter polytope} if all its dihedral angles are  submultiples 
of $\pi$. The group $G_F$ generated by  reflections in the facets of 
any Coxeter polytope $F$ is discrete, and $F$ is a fundamental domain 
of $G_F$.

On the other hand, any discrete group $G$ generated by reflections 
coincides with $G_F$ for some Coxeter polytope $F$. If $G=G_P$ for 
some non-Coxeter polytope $P$, then $P$ consists of several copies of 
$F$, and any two copies containing a facet in common are 
symmetric to each other with respect to this facet.\\  

\noindent
{\bf Spherical reflection groups}.
Let  $G$ be a reflection group acting on $\S^n$. 
Suppose that $G$ acts on $\S^n$ discretely, i.e. $G$ is a finite 
group.
Then  $G$ is generated by reflections 
in the facets of some spherical Coxeter polytope. 
It is shown by H.~S.~M.~Coxeter \cite{Coxeter} that 
any spherical Coxeter polytope containing no pair of antipodal points
of $\S^n$ is a simplex. 

To describe Coxeter polytopes we  use {\it Coxeter diagrams}.
A Coxeter diagram of a Coxeter polytope $F$
is a graph whose nodes $v_i$ correspond to the facets $\Pi_i$ of $F$.
Nodes $v_i$ and $v_j$ are joined by a $(k-2)$-fold edge 
if the dihedral angle
formed up by $\Pi_i$ and $\Pi_j$ equals $\frac{\pi}{k}$
(if $\Pi_i$ is orthogonal to $\Pi_j$, $v_i$ and $v_j$ are disjoint).
Indecomposable spherical and Euclidean Coxeter simplices
were classified by Coxeter~\cite{Coxeter}. The list of their 
Coxeter diagrams is represented  in Table~\ref{cox}.
\\




\noindent
{\bf Euclidean reflection groups}. 
Let  $G$ be a discrete reflection group acting on $\E^n$.
Let $F$ be a fundamental chamber of $G$.
As it is shown in~\cite{Coxeter}, $F$ is a direct product of several
simplices and simplicial cones. 

Suppose that $G$ is generated by a simplex (the simplex may
not be a Coxeter simplex). Then $G$ is an indecomposable group
and the fundamental chamber of $G$ is compact.
In this case the fundamental chamber of $G$ is one of the simplices 
 $\widetilde A_n$, $\widetilde B_n$, $\widetilde C_n$, $\widetilde D_n$,
$\widetilde E_6$, $\widetilde E_7$, $\widetilde E_8$, $\widetilde F_4$
and  $\widetilde G_2$ (see Table~\ref{cox}).

See~\cite{29}
for more information about discrete reflection groups.

We use the notation $A_n$, $\widetilde A_n$, $B_n$ and so on for Coxeter
simplices as well as for the groups generated by these simplices.
Finite root systems are denoted by $\mathrm A_n$,
$\mathrm B_n$ and so on. See~\cite{Kac} for background on root systems.    

Let $P$ be a simplex generating discrete reflection group
($P$ may not be a Coxeter simplex).
Clearly, in this case all the dihedral angles of $P$ are of the type $\frac{\pi m}{k}$.
Hence, for any simplex generating discrete reflection group
we can construct the following {\it generalized Coxeter diagram}:
the nodes $v_i$ of the diagram correspond to the facets $\Pi_i$ of $P$;
the nodes $v_i$ and $v_j$ are joined by a $(k-2)$-fold edge that is 
decomposed into $m$ parts
if the dihedral angle
formed up by $\Pi_i$ and $\Pi_j$ equals to $\frac{\pi m}{k}$.

\input{cox.txt}

\section{Families of simplices}

\noindent
{\bf Spherical simplices generating discrete reflection groups}.
Let $\E^{n}$ be the $n$-dimensional Euclidean  space,
and let $S^{n-1}$ be the unit sphere centered at the origin.
Any hyperplane of  $S^{n-1}$ is a section of $S^{n-1}$
by some  hyperplane of $\E^{n}$ containing the origin.
Hence, any $(n-1)$-dimensional simplex in  $S^{n-1}$ is an intersection 
of $S^{n-1}$ with an interior of some 
cone with $n$ facets in $\E^{n}$ centered at the origin. 

So, suppose that $P$ is a spherical $(n-1)$-dimensional simplex. 
Then we can define $P$ by unit outward normal vectors $f_1,...,f_n$ 
to the facets of the corresponding cone. 

Denote by $\Pi_1,...,\Pi_{n}$ the facets of $P$.
The hyperplanes containing $\Pi_1,...,\Pi_{n}$ decompose
 $S^{n-1}$ into $2^{n}$ simplices $P_1$, ..., $P_{2^n}$
encoded by $n$-tuples of vectors $\{\pm f_1,...,\pm f_{n}\}$. 
In this paper the set of simplices $P_1$,...,$P_{2^{n}}$
is called  a {\it family}.
Each of the simplices $P_1$,...,$P_{2^{n}}$ 
generates the same reflection group as $P$ does.
Thus, we can study the families instead of studying of simplices
 themselves. (In fact, any family contains at most $2^{n-1}$ simplices 
up to isometry: the simplex $\{f_1,...,f_{n}\}$ is always congruent to 
 $\{-f_1,...,-f_{n}\}$).

Let $P$ be a simplex generating a discrete reflection group $G_P$. 
Clearly, the dihedral angles of $P$ are rational 
numbers multiplied by $\pi$. Moreover, if $G$ is an indecomposable
spherical reflection group different from  $G_2^{(m)}$, $H_3$ and $H_4$, 
then any dihedral angle of $P$
is either $\frac{\pi}{k}$ or  $\frac{\pi(k-1)}{k}$, where $k=2,3$ or $4$ 
(cf. Table~\ref{cox}).

Analogously to Coxeter simplices, any simplex  
whose dihedral angles equal  either 
 $\frac{\pi}{k}$ or  $\frac{\pi(k-1)}{k}$
can be represented by the following {\it family diagram}:
the nodes $v_i$ correspond to the facets of $P$;
the nodes $v_i$ and $v_j$ are joined by a $(k-2)$-fold edge 
if the angle
formed up by $\Pi_i$ and $\Pi_j$ is either $\frac{\pi}{k}$
or  $\frac{\pi(k-1)}{k}$
(if $\Pi_i$ is orthogonal to $\Pi_j$, $v_i$ and $v_j$ are disjoint).

Note, that if $P$ is a Coxeter simplex then the  diagram corresponding
 to $P$ is a Coxeter diagram of $P$.

Any two simplices in one family
have the same family diagram. 
So, we can assign to a family the diagram of any simplex contained 
in the family. We call this diagram a {\it family
diagram}. It is proved in~\cite{sph} (Lemma 2) that any graph having 
no edge of
multiplicity greater than two is a family diagram for at most one
family of simplices. The same statement for the dihedral group 
$G_2=G_2^{(6)}$ is evident.

\vspace{8pt}
\noindent
{\bf Euclidean simplices generating discrete reflection groups}.
Let $P$ be a simplex in $\E^n$
generating discrete reflection group $G_P$.
Let $\Pi_0,...,\Pi_{n}$ be the facets of $P$,
and $f_0,...,f_n$ be the outward unit normal vectors to 
$\Pi_0,...,\Pi_{n}$.
The hyperplanes containing $\Pi_0,...,\Pi_{n}$ decompose
$\E^{n}$ into $2^{n+1}-1$ domains $P_1$,...,$P_{2^{n+1}-1}$
encoded by the vectors $\{\pm f_0,...,\pm f_{n}\}$. 
All these domains except the initial simplex $P$ are non-compact.
Simplex $P$ is encoded by the vectors $\{f_0,..., f_{n}\}$.
Each  of the domains $P_i$ generates the same reflection group as $P$ does.

By a {\it family of Euclidean simplices} we call a set of simplices
having the same set of vectors $\{\pm f_0,...,\pm f_{n}\}$.
Note, that in the Euclidean case any two simplices contained in one
family are mutually similar (moreover, any two of these simplices are 
homothetic).

We classify Euclidean simplices generating discrete reflection
groups up to similarity.

Furthermore, note that any dihedral angle of $P$ is either
$\frac{\pi}{k}$ or  $\frac{\pi(k-1)}{k}$,
where $k=2,3,4$ or $6$ (see Table~\ref{cox}).
Hence, we can define {\it family diagrams} for Euclidean simplices 
in the same way as for  the spherical ones.

\begin{lemma}
\label{unique}
Let $\Phi$ and $\Psi$ be two different families of Euclidean simplices
generating discrete reflection groups. Then their family 
diagrams are distinct.

\end{lemma}

The proof of the lemma follows the proof of Lemma~2 from
paper~\cite{sph}. 

\section{Special vertices}
In this section we prove several auxiliary facts.

A hyperplane $\alpha$ is called a {\it mirror} of the group $G$  
if $G$ contains a reflection with respect to $\alpha$.

\begin{lemma}
\label{parallel}
Let $P$ be a simplex in $\E^n$ generating discrete reflection
group $G_P$.
Then no mirror of $G_P$ decomposing $P$ is
 parallel to a facet of $P$.
 
\end{lemma}

\begin{proof}
Let $\Pi_0,...,\Pi_n$ be the facets of $P$ and let $V_0$ be the vertex 
opposite to $\Pi_0$.
Suppose that there exists a mirror of $G_P$ that is parallel to $\Pi_0$ and
decomposes $P$. Since $G_P$ is discrete, there exist only finitely
many of such mirrors. Let $m$ be one of them closest to $V_0$.

Denote by $h$ the homothety with center $V_0$ taking  $\Pi_0$ to $m$.
Denote by $r_i$ the reflection with respect to $\Pi_i$ ($i=0,...,n$), 
and by $r$ the reflection with respect to $m$.
Since $r\in G_P$,  $r=r_{i_1}...r_{i_l}$ for some $l$.
From the other hand, $r=hr_0h^{-1}$.

Consider the reflection
$hrh^{-1}=(hr_{i_1}h^{-1})...(hr_{i_l}h^{-1})$.
Since $hr_{i}h^{-1}\in G_P$ for any $i=0,...,n$,  
we have $hrh^{-1}\in G_P$. Furthermore,  $hrh^{-1}$ is a reflection
with respect to some hyperplane $m'$ parallel to $m$.
Moreover, $m'$ decomposes $P$ and goes closer to $V_0$ than $m$. 
This contradicts to the choice of $m$.

\end{proof}

\begin{lemma}
\label{Fix}
Let $P$ be a simplex generating  discrete group $G_P$.
Let $V$ be a vertex of $P$ and  $\Pi_1,...,\Pi_n$ be the facets of
$P$ containing $V$.
Then the stabilizer $Fix(V,G_P)$ of $V$ in  $G_P$ coincides with
the group $G$
generated by the reflections with respect to $\Pi_1,...,\Pi_n$.
 
\end{lemma}

\begin{proof}

Suppose that the lemma is false. 
Then there exists a non-empty set $M$ of  simplices 
for which the statement of the lemma is broken.
We may assume that $P\in M$ is a simplex minimal by the inclusion.
Let $V$ be a vertex of $P$ for which the statement of the lemma is false.
Then among the dihedral angles formed up by $\Pi_1,...,\Pi_n$ 
 there exists a dihedral angle $\pi \frac{k}{m}$, where $1<k<m$
are mutually co-prime integers. Suppose that this angle is 
formed up by $\Pi_1$ and $\Pi_n$.

Consider a mirror $\Pi$ of $G$ such that $\Pi$ contains 
$\Pi_1\cap \Pi_n$, and
the angle formed up by $\Pi$ and $\Pi_1$ is equal to $\frac{\pi}{m}$.
This mirror decomposes $P$ into two simplices $P_1$ and $P_2$.
Let $P_1$ be the simplex having a facet $\Pi_1$. 
Denote by $G_{P_1}$ the group generated by $P_1$. Clearly, $G_{P_1}=G_P$.
Furthermore, the group generated by the reflections with respect to the
facets $\Pi,\Pi_1,...,\Pi_{n-1}$ coincides with 
the group generated by the reflections with respect to the
facets $\Pi_1,...,\Pi_n$. From the other hand, the stabilizer
 $Fix(V,G_{P_1})$ of $V$ in $G_{P_1}$ coincides with $Fix(V,G_{P})$.
By the assumption, $Fix(V,G_{P})$ does not coincide with the group generated by
the reflections with respect to $\Pi_1,...,\Pi_n$.
Hence,   $Fix(V,G_{P_1})$ differs from the group generated by the
reflections with respect to $\Pi,\Pi_1,...,\Pi_{n-1}$.
Thus, the statement of the lemma is broken for $P_1$.
This contradicts to the assumption that $P\in M$ is the minimal
simplex, and the lemma is proved.

\end{proof}

Let $P$ be any Euclidean polytope generating discrete reflection group. 
Suppose that
there exists a vertex $V$ of $P$ such that the stabilizer of
$V$ contains a linear part of any element of $G_P$.
We call such a vertex a {\it special vertex } of $P$.
It is known  (see~\cite{Kac}, Ch.~6) that any Euclidean Coxeter
simplex has at least one special vertex. 

\begin{lemma}
\label{spec}
Let $P$ be a simplex generating  discrete reflection group $G_P$. 
Then $P$ has at list one special vertex.

\end{lemma}

\begin{proof}
Let $F$ be a fundamental polytope of $G_P$ contained in $P$.
Since $G_P$ is indecomposable, $F$ is a simplex. Any fundamental
polytope of discrete reflection group is Coxeter polytope, thus $F$ is
a Coxeter simplex. Let $V$ be a special vertex of $F$.

Suppose that $F$ is either an inner point of $P$ or an inner point of
some face of $P$. Let $\Pi$ be any facet of $P$ not containing $V$.
Let $m$ be the mirror of $G_P$ parallel to $\Pi$ and containing $V$.
Then $m$ cuts $P$, that contradicts to Lemma~\ref{parallel}.

Thus, $V$ is a vertex of $P$. Let $\Pi_1,...,\Pi_n$ be the facets of $P$
containing $V$.  
By Lemma~\ref{Fix}, the reflections with respect to  $\Pi_1,...,\Pi_n$
generate the stabilizer $Fix(V,G_P)$ of $V$ in $G_P$.
Thus, $V$ is a special vertex of $P$.

\end{proof}

\section{Simplices generating given reflection group}
\label{algorithm}

Let $G$ be an indecomposable Euclidean reflection group.  
Let $P$ be a simplex generating $G$, and $f_0,f_1,...,f_n$ 
be the vectors orthogonal to the facets of $P$.
By Lemma~\ref{spec}, there exists a special vertex $V$ of $P$.
Let  $\Pi_0,...,\Pi_n$ be facets of $P$, such that $\Pi_i$ is orthogonal
 to $f_i$ for $i=0,\dots,n$. We may assume that $V$ is contained in  
$\Pi_1,...,\Pi_n$. 

The stabilizer $Fix(V,G)$ of $V$ in $G$ is a Weyl 
group of some finite root system $\Delta$. Suitably normalizing the 
vectors $f_0,f_1,...,f_n$,
we may assume that these vectors belong to $\Delta$. Since $V$ is a 
special vertex of $P$, there exists a mirror of $G$ through $V$ 
parallel to $\Pi_0$. By Lemma~\ref{Fix}, suitably normalizing $f_0$ 
we may assume that $f_0$ is contained in $\Delta$.    
Note, that $f_0,f_1,...,f_n$ are linearly dependent vectors,
however any $n$ of these vectors are linearly independent.

Now we are able to find all families of simplices generating 
given group $G$.

Let $F$ be a fundamental simplex of $G$,
and let $U$ be a special vertex of $F$.
Denote by $W$ the stabilizer $Fix(U,G)$, and let $\Delta$ be the 
corresponding root system.    
Let $v_0,...,v_n$ be any set containing $n+1$ vectors of
$\Delta$, such that any $n$ of these vectors are linearly independent.
Let  $\Pi_1,...,\Pi_n$ be the mirrors of $G$ containing $U$ and
orthogonal to $v_1,...,v_n$.
Let $\Pi_0$ be a mirror of $G$ closest to $U$, orthogonal to $v_0$ and
not containing  $U$.
Then the mirrors $\Pi_0,\Pi_1,...,\Pi_n$ are the facets of some 
simplex $P$
generating a finite index subgroup of $G$.
If, in addition, the reflections with respect to 
 $\Pi_1,...,\Pi_n$ generate $W$, 
then some simplex similar to $P$ generates $G$ 
(the only exception is the group $W=B_n=C_n$; in this case
we can obtain either $\widetilde B_n$ or $\widetilde C_n$
depending on the vector $v_0$, see Section~\ref{SectionBC}).

In other word, to classify families generating $G$ 
it is sufficient to follow the algorithm:

\begin{itemize}
\item[1)] Find all linearly independent systems
$f_1,...,f_n$ in $\Delta$, generating $W$
(we say that $f_1,...,f_n$ generate  $W$,
if the reflections with respect to the hyperplanes through the origin
orthogonal  to these vectors generate $W$).

\item[2)]
For each system $f_1,...,f_n$ obtained on the previous step, add a vector
$f_0\in \Delta$ in such a way that any $n$ of vectors 
 $f_0,f_1,...,f_n$ are linearly independent.
This vector $f_0$ should be chosen by all possible ways.

\item[3)] If $W=B_n=C_n$, examine which of the groups
$\widetilde B_n$ and  $\widetilde C_n$ is generated by
$f_0,f_1,...,f_n$. 

\item[4)] Among the obtained systems  
one should find the systems corresponding to different families,
i.e.  families having different family diagrams.

\end{itemize}

The order of vectors in the systems is not important for us.
The Weyl group acts on  $\E^n$,
hence, it acts on $(n+1)$-tuples of vectors.
We do not differ   $(n+1)$-tuples equivalent with respect to
this action. 

\subsection{Simplices generating $\widetilde A_n$}
\label{SecA}

Let $P$ be a simplex generating the group $\widetilde A_n$.
We may assume that the vectors $f_0,f_1,...,f_n$ belong to the root 
system 
$\mathrm A_n=\{\pm(h_i-h_j)\}, 0\le i<j\le n$, where
$h_0,...,h_n$ is a standard basis of $\E^{n+1}$).

For any simplex $P=\{f_0,f_1,...f_n\}$ generating $\widetilde A_n$
we construct the following {\it graph} $\,\Gamma(P)$:
the nodes $v_0,...,v_n$ of  $\Sigma$ correspond to the vectors 
$h_0,...h_n$; the nodes $v_i$ and $v_j$ are joined by an edge
if one of vectors $(h_i-h_j)$ and $-(h_i-h_j)$ belongs to the set 
$\{f_0,...,f_n\}$. 
Clearly, two simplices from one family have the same 
graph. 

\begin{theorem}
\label{sectionA}
There exists a unique family of simplices generating 
$\widetilde A_n$.
This  family consists of Coxeter simplices $\widetilde A_n$.
\end{theorem}

\begin{proof}
Let $P=\{f_0,f_1,...f_n\}$ be a simplex generating  $\widetilde A_n$
and  $\Gamma(P)$ be the corresponding graph.
Then $\Gamma(P)$ contains exactly $n+1$ nodes and the same number
of edges.
Since any $n$ vectors contained in the set  $\{f_0,...,f_n\}$ are
linearly independent, 
 $\Gamma(P)$ has no cycles containing less than $n+1$ nodes.
Hence,  $\Gamma(P)$ is a cycle with $n+1$ nodes.

Note, that the graph  $\Gamma(P)$ determines the family diagram 
of  $P$. In more details, the nodes of the family diagram
correspond to the edges of   
$\Gamma(P)$, two nodes are joined if the corresponding edges
of  $\Gamma(P)$  have a common node.
Hence, in case of $\widetilde A_n$ the family diagram is a cycle,
too. By Lemma~\ref{unique},   
the family of simplices is completely determined by a family diagram.
Thus, $P$ is a Coxeter simplex  $\widetilde A_n$.

\end{proof}

\subsection{Simplices generating $\widetilde B_n$ and $\widetilde C_n$} 
\label{SectionBC}

Let $P$ be a simplex generating the group $\widetilde B_n$ or 
$\widetilde C_n$.
We may assume that the vectors $f_0,f_1,...,f_n$ belong to the root 
system 
$\mathrm B_n=\{\pm h_i,\ \pm h_i \pm h_j\}, 1\le i<j\le n$, where
$h_1,...,h_n$ is a standard basis of $\E^{n}$).

For any simplex $P=\{f_0,f_1,...f_n\}$ generating $\widetilde B_n$ or 
$\widetilde C_n$ 
we construct the following {\it graph} $\,\Gamma(P)$:
the nodes $v_1,...,v_n$ of  $\Sigma$ correspond to the vectors 
$h_1,...h_n$; the nodes $v_i$ and $v_j$ are joined by an edge
if one of  vectors $\pm (h_i-h_j)$ and $\pm (h_i+h_j)$ belongs to 
the set 
$\{f_0,...,f_n\}$; the node $v_i$ is {\it marked} if   
$\{f_0,f_1,...,f_n\}$ contains  $\pm h_i$.
If $\{f_0,f_1,...,f_n\}$ contains both
 $\pm(h_i+h_j)$ and $\pm(h_i-h_j)$,   
the nodes $v_i$ and $v_j$ are joined by two edges. 

Since the system of vectors $f_0,f_1,...,f_n$  is indecomposable, 
 $\Gamma(P)$ is connected.
Clearly, $\Gamma(P)$ contains at least one marked node
(otherwise $f_i$ belongs to $\mathrm D_n$ for any $i$,
where  $\mathrm D_n$ is embedded in $\mathrm B_n$ 
as a set of long roots).

\begin{lemma}
\label{marked}
Let $P$ be a simplex generating $\widetilde B_n$ or 
$\widetilde C_n$. \\
1) If $\Gamma(P)$ contains more than one marked node 
then $P$ is a Coxeter simplex $\widetilde C_n$.\\
2) If $\Gamma(P)$ contains a unique
marked node then $\Gamma(P)$ is one of the  graphs shown in 
the left column of Table~\ref{beat}.

\end{lemma}

\begin{figure}[!h]
\begin{center}
\epsfig{file=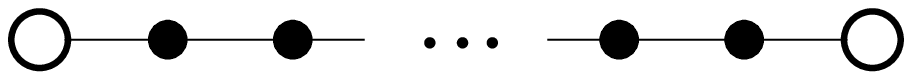,width=0.4\linewidth}
\vspace{-10pt}
\end{center}
\caption{Graph $\Gamma(\widetilde C_n)$.}
\label{Sigma(C_n)}
\end{figure}

\begin{proof}
Suppose that $\Gamma(P)$ contains more than one marked node.
Then $\Gamma(P)$ contains a path from one marked node to another.
Consider such a path that does not intersect itself.
The vectors corresponding to the edges and marked nodes of this
path are linearly dependent. Hence, the path contains all edges of
 $\Gamma(P)$, and  $\Gamma(P)$ is the graph shown in
Fig.~\ref{Sigma(C_n)}.
In  this case $P=\widetilde C_n$. 

Now suppose that $\Gamma(P)$ contains a unique marked node.
Since the number of edges of $\Gamma(P)$ equals $n$
and the number of nodes equals $(n+1)-1=n$,  
$\Gamma(P)$ contains a unique cycle
(it may consist of two nodes).
Since any $n$ vectors are linearly independent, if the cycle contains
less than $n$ edges then it does not contain the marked node.
Hence, $\Gamma(P)$ have a subgraph shown in the left column of
Table~\ref{beat}.
From the other hand, this subgraph corresponds to linearly
dependent system of vectors.
Hence, this subgraph coincides with $\Gamma(P)$.

\end{proof}

\begin{theorem}
\label{B}
There exists exactly $n-1$ families of simplices generating
the group $\widetilde B_n$. Family diagrams and  generalized 
Coxeter diagrams of these simplices are shown in Table~\ref{beat}.

\end{theorem}

\begin{proof}
By the second statement of Lemma~\ref{marked},
any simplex $P$ generating $\widetilde B_n$ 
corresponds to a graph $\Gamma(P)$ shown in the left column of
Table~\ref{beat}. 

Let us show that any graph $\Gamma$  shown in the left column of
Table~\ref{beat} corresponds to a simplex generating $\widetilde B_n$.
It is easy to see that any $n$ of vectors  $f_0,f_1,...,f_n$ 
are linearly independent and all these vectors are linearly dependent.
Hence,  $\Gamma$ corresponds to some simplex $P$ in $\E^n$.
The group generated by $P$ is a maximal rank indecomposable subgroup of
$\widetilde B_n$, and the fundamental simplex of this subgroup has a
dihedral angle equal to $\frac{\pi}{4}$. 
By~\cite{eucl}, $P$ generates either $\widetilde B_n$ or 
$\widetilde C_n$. 

Let $\varphi$ and $\psi$ be two dihedral angles equal to $\frac{\pi}{4}$
formed up by mirrors of $G_P$.
Then there exists an element $\gamma$ of $G_P$ such that 
$\gamma(\varphi)=\psi$.
The group  $\widetilde C_n$ contains two equivalency classes of such
dihedral angles. Hence, $G_P\ne \widetilde C_n$ and
$P$ generates $\widetilde B_n$.

Furthermore, let us show that each graph shown in the left column
of the Table~\ref{beat} corresponds to a unique family of simplices
generating $\widetilde B_n$.
Indeed, the family diagram of the family containing $P$ can be
easily recovered from  $\Gamma(P)$:
all but one nodes of the diagram correspond to the edges  
of $\Gamma(P)$,
two nodes are adjacent if the corresponding edges 
have a common point;
the rest node corresponds to the marked node  
of $\Gamma(P)$, this node is joined  by a 2-fold
edge with all the nodes that correspond to edges of $\Gamma(P)$
incident to the marked node. 
Thus, any graph shown in the left column of Table~\ref{beat}
corresponds to a unique family diagram, and, by Lemma~\ref{unique}, 
it corresponds to a unique family.
An explicit calculation shows that the generalized Coxeter diagram
of simplex generating $\widetilde B_n$ is one of the diagrams shown
in the right column of Table~\ref{beat}.

Thus,  families of  simplices generating $\widetilde B_n$
are in one-to-one correspondence with  graphs shown in the left
column of Table~\ref{beat}.
To find the number of these families note, that 
$\Gamma(P)$ is uniquely determined by the number of edges in the
cycle. The latter is any integer number $N$ satisfying $2\le N \le n$.

\end{proof}

\begin{table}[!h]
\caption{Simplices generating $\widetilde B_n$. }
\label{beat}
\vspace{5pt}
\begin{center}
\begin{tabular}{|c|c|c|}
\hline
\ $\Gamma(P)$\ &\  Family diagram\ 
& 
\begin{tabular}{c}
\raisebox{0pt}[15pt][10pt]{Generalized} \\
\raisebox{5pt}[5pt][5pt]{Coxeter diagram}
\end{tabular}
\\
\hline
&&\\
\epsfig{file=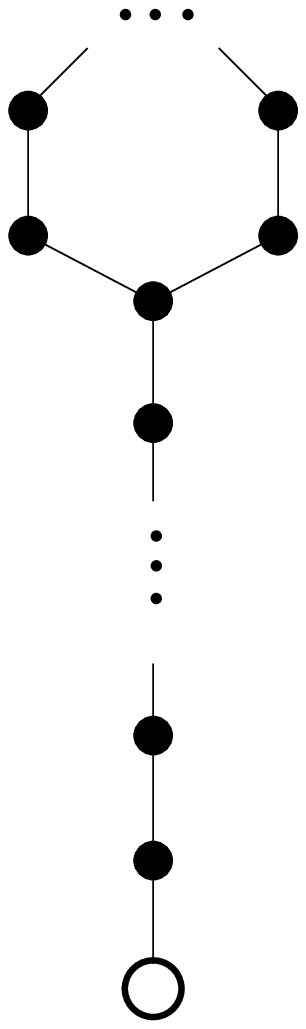,width=0.092\linewidth}&
\epsfig{file=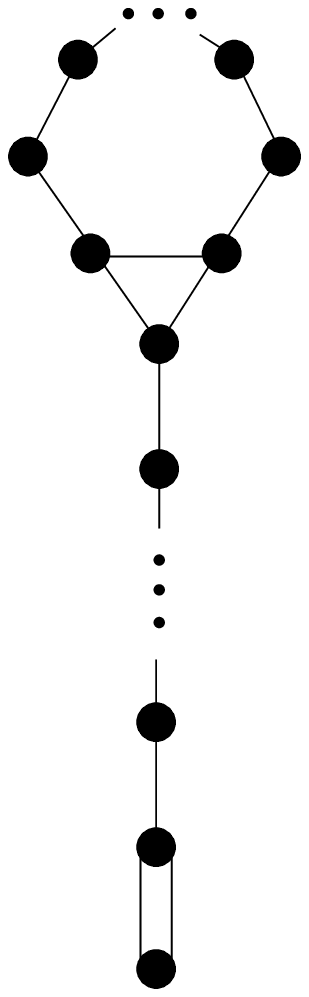,width=0.098\linewidth}&
\epsfig{file=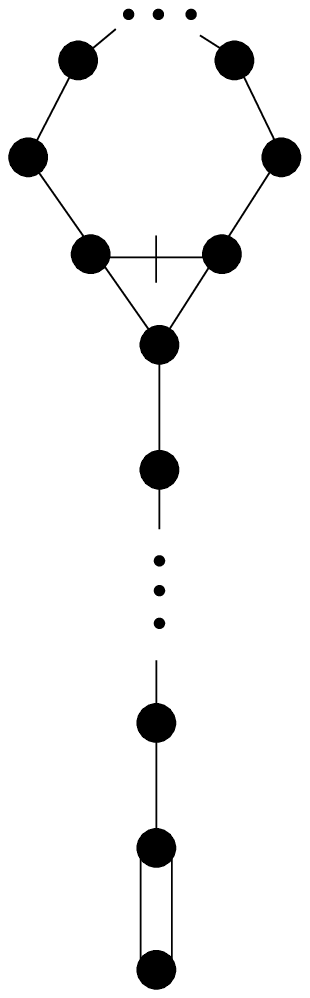,width=0.098\linewidth}\\
&&\\
\hline
&&\\
\raisebox{10pt}[0pt][0pt]{\epsfig{file=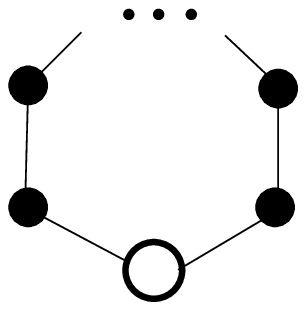,width=0.101\linewidth}}&
\raisebox{6pt}{\epsfig{file=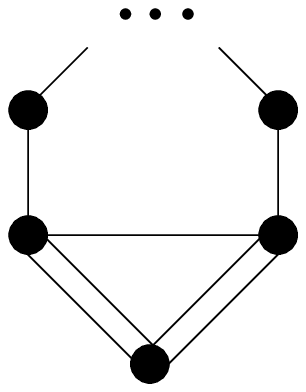,width=0.095\linewidth}}&
\raisebox{6pt}{\epsfig{file=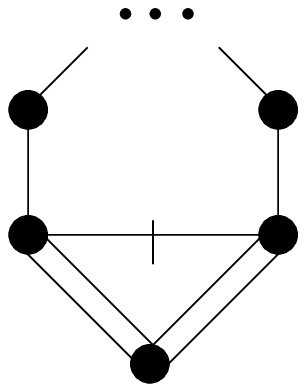,width=0.095\linewidth}}\\
\hline
\end{tabular}
\end{center}
\end{table}

\noindent
As a corollary of Theorem~\ref{B} and the first statement of 
Lemma~\ref{marked} we obtain the following theorem:

\begin{theorem}
\label{C}
There exists a unique family of simplices generating 
$\widetilde C_n$.
This  family consists of Coxeter simplices $\widetilde C_n$.

\end{theorem}

\subsection{Simplices generating $\widetilde D_n$}
\label{sectionD}

Let $P$ be a simplex generating the group $\widetilde D_n$.
We may assume that the vectors $f_0,f_1,...,f_n$ belong to the 
root system 
$\mathrm D_n=\{\pm h_i \pm h_j)\}, 1\le i<j\le n$, where
$h_1,...,h_n$ is a standard basis of $\E^{n}$).

For any simplex $P=\{f_0,f_1,...f_n\}$ generating $\widetilde D_n$
we construct the following {\it graph} $\,\Gamma(P)$:
the nodes $v_1,...,v_n$ of  $\Sigma$ correspond to the vectors 
$h_1,...h_n$; the nodes $v_i$ and $v_j$ are joined by an edge
if one of vectors $\pm (h_i-h_j)$ and $\pm (h_i+h_j)$ belongs to the set 
$\{f_0,...,f_n\}$; if $\{f_0,f_1,...,f_n\}$ contains both $\pm(h_i+h_j)$ and  
$\pm(h_i-h_j)$ then $v_i$ and $v_j$ are joined by two edges.

The system of vectors $f_0,f_1,...,f_n$  is indecomposable,
hence, $\Gamma(P)$ is connected.
Since  $\Gamma(P)$ has $n$ nodes and $n+1$ edges,   
 $\Gamma(P)$ contains at least two cycles $C_1$ and $C_2$.

Suppose that  $C_1$ and $C_2$ have no common nodes.
Since the graph containing two disjoint cycles corresponds to 
a linearly dependent system of vectors, 
$\Gamma(P)$ is one of the graphs shown in the left column 
of Table~\ref{veslo}. 
If $C_1$ and $C_2$ have a unique common node then
$\Gamma(P)$ is the graph shown at the bottom of left column
of Table~\ref{veslo}. 

\begin{lemma}
\label{2cycles}
Suppose that $\Gamma(P)$ contains two cycles having at least two
common nodes. Then the system of vectors  $f_0,f_1,...,f_n$
contains $n$ linearly dependent vectors.

\end{lemma}

\begin{proof}
Consider the graph $\Gamma(P)$ colored by the following way: edges
corresponding to vectors 
$\pm(h_i+h_j)$ are red, and the rest edges, i.e. the edges
corresponding to $\pm(h_i-h_j)$, are black. 
Note that substituting the vector $h_i$  by $-h_i$ we change the 
color of all edges incident to  $v_i$.
Thus, preserving the vectors  $f_0,f_1,...,f_n$,
we can make all but one edge of a given cycle black
(the rest edge is either red or black). 
Vectors $h_1-h_2,h_2-h_3,...,h_k-h_1$ are linearly dependent.
Hence, each cycle containing an even number of  red edges
corresponds to a system of linearly dependent vectors.

Consider common nodes of two cycles contained in $\Gamma(P)$.
There are at least three paths  $L_1,L_2$ and $L_3$ joining these
nodes in $\Gamma(P)$.
Denote by $c(L_i)$ the number of red edges in $L_i$.
Then $c(L_i)+c(L_j)$ is even for some $i\ne j$, $i,j=1,2,3$.
We assume that $c(L_1)+c(L_2)$ is even. Then the cycle $C=L_1\cup L_2$
contains even number of  red edges, so, it corresponds
to some linearly dependent vectors. 
Since some edges of $\Gamma(P)$ do not belong to $C$,
the number of these linearly dependent vectors is less than $n+1$.
The contradiction  proves the lemma.

\end{proof}

\begin{cor}
\label{d}
$\Gamma(P)$ contains exactly two cycles and coincides with one of 
graphs shown in the left column of Table~\ref{veslo}.

\end{cor}

Thus, we obtain the following theorem:

\begin{theorem}
The group $\widetilde D_n$ is generated by exactly  
 $\frac{1}{4}n(n-2)$ families of simplices if $n$ is even,
and by exactly $\frac{1}{4}(n-1)^2$ families if $n$ is odd.
Family diagrams and generalized Coxeter diagrams 
of these simplices are presented in Table~\ref{veslo}.

\end{theorem}

\begin{table}[!h]
\caption{Simplices generating $\widetilde D_n$. }
\label{veslo}
\vspace{5pt}
\begin{center}
\begin{tabular}{|c|c|c|}
\hline
\ $\Gamma(P)$\ &\  Family diagram\ 
& 
\begin{tabular}{c}
\raisebox{0pt}[15pt][10pt]{Generalized} \\
\raisebox{5pt}[5pt][5pt]{Coxeter diagram}
\end{tabular}
\\
\hline
&&\\
\epsfig{file=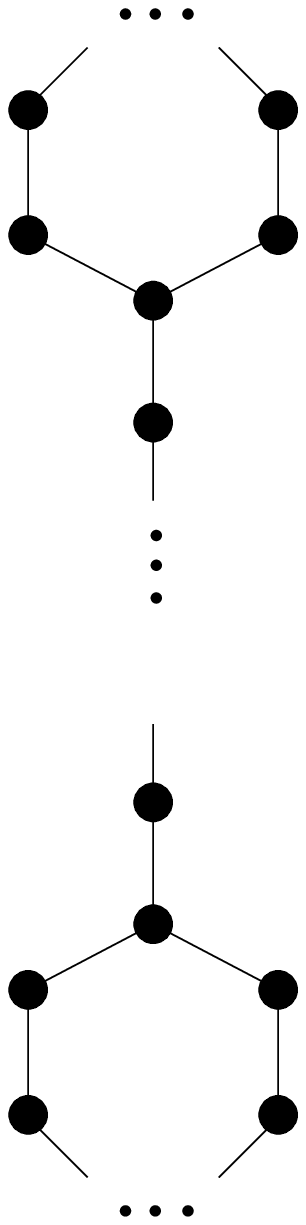,width=0.0939\linewidth}&
\epsfig{file=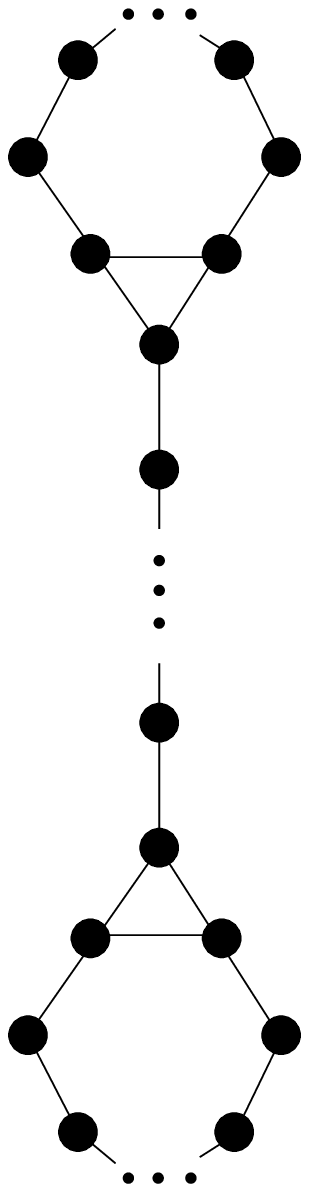,width=0.098\linewidth}&
\epsfig{file=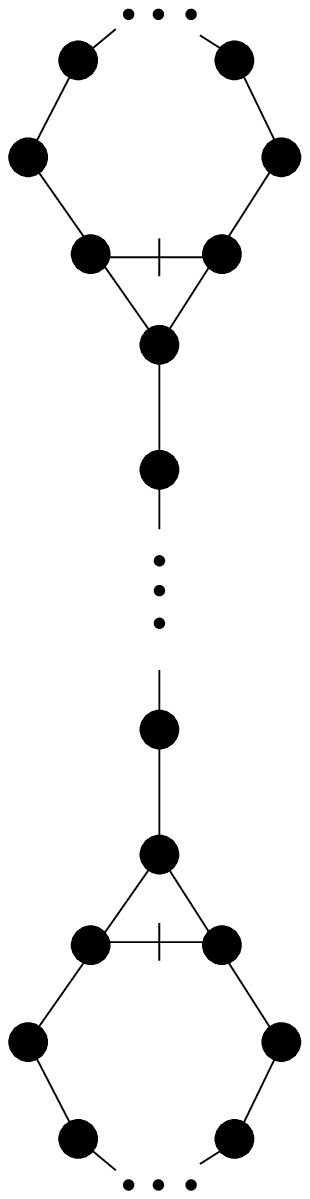,width=0.098\linewidth}\\
&&\\
\hline
&&\\
\raisebox{7pt}[0pt][0pt]{\epsfig{file=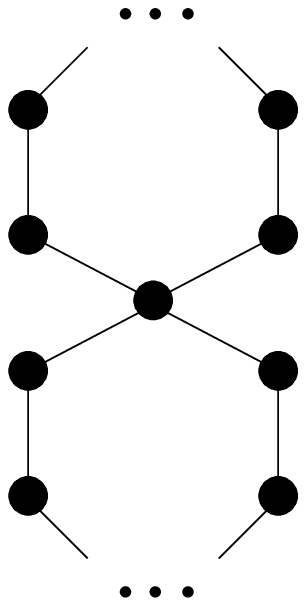,width=0.097\linewidth}}&
\raisebox{6pt}{\epsfig{file=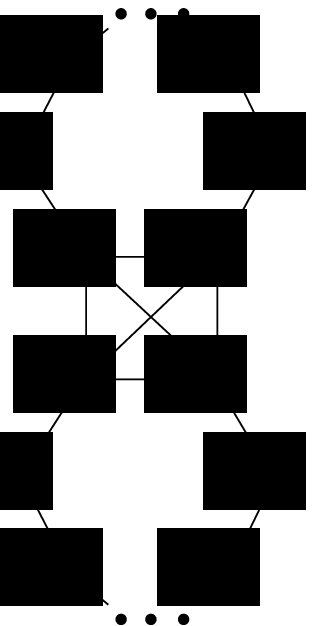,width=0.095\linewidth}}&
\raisebox{6pt}{\epsfig{file=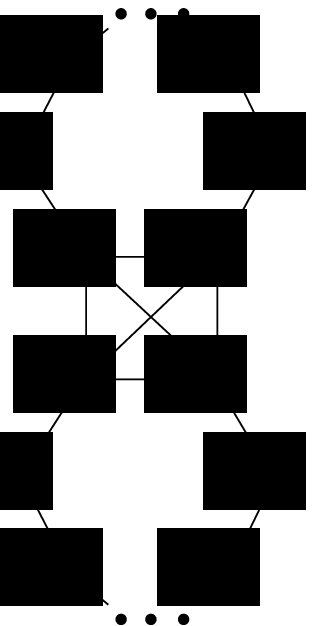,width=0.095\linewidth}}\\
\hline
\end{tabular}
\end{center}
\end{table}

\begin{proof}

By Cor.~\ref{d}, a simplex $P$ generating $\widetilde D_n$ 
corresponds to a graph $\Gamma(P)$ shown in the left column of
Table~\ref{veslo}. 

Let us show that any graph $\Gamma$  shown in the left column of
Table~\ref{veslo} corresponds to a simplex generating $\widetilde
D_n$.
For each cycle of $\Gamma$ choose an edge $v_iv_j$ and put the vector
 $h_i-h_j$ in correspondence with this edge (in other words, suppose
these edges to be red). For all other edges $v_kv_l$  take  vectors $h_k+h_l$.
It is easy to see that any $n$ of these vectors   
are linearly independent and all these vectors are linearly dependent.
Hence,  $\Gamma$ corresponds to some simplex $P$ in $\E^n$.
By~\cite{eucl}, the group  $\widetilde D_n$ has no finite index
indecomposable subgroups different from $\widetilde D_n$. 
Therefore, $P$ generates  $\widetilde D_n$. 

Now, show that each graph $\Gamma$ shown in the left column of
Table~\ref{veslo} corresponds to a unique family of simplices generating
 $\widetilde D_n$.
Indeed, a family diagram  of the family containing $P$ 
can be recovered from $\Gamma(P)$: nodes of diagram 
correspond to  edges of  $\Gamma$, two nodes are adjacent if the
corresponding edges have a common point.
Therefore,  the family does not depend on the choice of initial edges. 
Moreover, in the beginning of the procedure we could make red not one but
several edges: the family diagram would not be changed.  
Hence, family diagram of simplices generating $\widetilde D_n$
are in one-to-one correspondence with  graphs shown in the left
column of Table~\ref{veslo}.
By Lemma~\ref{unique}, each family diagram corresponds to a unique 
family.
Thus, families of simplices generating  $\widetilde D_n$
are in one-to-one correspondence with graphs shown in the left
column of Table~\ref{veslo}. 
An explicit calculation shows that generalized Coxeter diagram
of simplex generating $\widetilde D_n$ is one of diagrams shown
in the right column of Table~\ref{veslo}.

To find the number of families note, that 
$\Gamma(P)$ is uniquely determined by the numbers $N_1$ and $N_2$ of
edges in two cycles. 
Numbers $N_1$ and $N_2$ are 
any two integers satisfying $N_1+N_2\le n+1$ and
$2\le N_1\le N_2 \le n$. The number of such pairs $(N_1,N_2)$
equals either $\frac{1}{4}n(n-2)$  
or $\frac{1}{4}(n-1)^2$ if $n$ is even or odd respectively.

\end{proof}

\subsection{Simplices generating other groups}

In Sections~\ref{SecA}--\ref{sectionD} we have described all 
families generating the groups
$\widetilde A_n$, $\widetilde B_n$,
$\widetilde C_n$  and $\widetilde D_n$.
Now, we are left to classify  families generating finitely many of
other indecomposable Euclidean reflection groups.
Namely, we a left with the groups 
$\widetilde E_6$, $\widetilde E_7$,
$\widetilde E_8$, $\widetilde F_4$,  
and $\widetilde G_2$. 
We can find the complete answer following the algorithm
contained in the beginning of Section~\ref{algorithm}.
As the result we obtain the lists which are rather large:
there exist \\
\phantom{aaaaaaaa} 17 families of simplices generating
  $\widetilde E_6$,\\ 
\phantom{aaaaaaaa} 142 families of simplices generating $\widetilde E_7$, \\
\phantom{aaaaaaaa} 1736 families of simplices generating $\widetilde E_8$,\\
\phantom{aaaaaaaa} 11 families of simplices generating $\widetilde F_4$, \\
\phantom{aaaaaaaa} and 2 families of simplices generating $\widetilde G_2$. \\
Appendix contains the complete list of these families.

\section*{Appendix}

Appendix contains the list of families generating
$\widetilde E_6$, $\widetilde E_7$,
$\widetilde E_8$, $\widetilde F_4$,  
and $\widetilde G_2$.
Families generating $E_6$, $E_7$ and $E_8$ are encoded
in the following way. For a family $\pm f_0,\pm f_1,...,\pm f_n$
construct a symmetrical matrix $G^+=\{g_{i,j}\}$, where
$g_{i,j}=2|(f_i,f_j)|$.
This is a doubled unsigned Gram matrix of the system
$f_0, f_1,..., f_n$.
The upper triangle of $G^+$ is filled up by $0$ and $1$.
Let $p$ be a decimal number which is equal to
the binary number
$\overline{
g_{1,2}g_{1,3}...g_{1,n}g_{2,3}...g_{2,n}...g_{i,i+1}...g_{i,n}...g_{n-1,n}}\!$.
The number $p$ depends on the ordering of vectors $f_0,f_1,...f_n$.
We choose $p$ the smallest possible.

\subsection*{Families generating $\widetilde E_6$}
\input{res_e6.txt}

\subsection*{Families generating $\widetilde E_7$}
\input{res_e7.txt}

\subsection*{Families generating $\widetilde E_8$}
\input{res_e8.txt}

\subsection*{Families generating $\widetilde F_4$}

 Families generating $F_4$ are encoded
in the following way. For a family determined by 
$\pm f_0,\pm f_1,...,\pm f_4$
construct a symmetric matrix $G^+=\{g_{i,j}\}$,
where $g_{i,j}=0$ if $f_i$ is orthogonal to $f_j$,
 $g_{i,j}=1$ if $\angle f_if_j=\frac{\pi}{3}$ or
$\frac{2\pi}{3}$,
 $g_{i,j}=2$ if $\angle f_if_j=\frac{\pi}{4}$ or
$\frac{3\pi}{4}$.
Let $p$ be a decimal number which is equal to
the base three number
$\overline{
g_{1,2}g_{1,3}...g_{1,n}g_{2,3}...g_{2,n}...g_{i,i+1}...g_{i,n}...g_{n-1,n}}$.
The number $p$ depends on the numbering of vectors $f_1,...f_n$.
We choose $p$ the smallest possible.
Then $p$ depends only on the family of simplices.

\input{res_f4.txt}

\subsection*{Families generating $\widetilde G_2$}

The group $G_2$ can be generated either by triangle with angles
$(\frac{\pi}{2},\frac{\pi}{3},\frac{\pi}{6})$ or
by triangle with angles
$(\frac{2\pi}{3},\frac{\pi}{6},\frac{\pi}{6})$.

\vspace{25pt}
\noindent
Independent University of Moscow, Russia \\
Max-Planck Institut f\"ur Mathematik Bonn, Germany\\
e-mail: \phantom{ow} felikson@mccme.ru\qquad pasha@mccme.ru

\end{document}

%% file: cox.txt

\begin{table}
\caption{{Coxeter diagrams.} Connected elliptic and parabolic Coxeter diagrams are
listed in left and right columns respectively.
 Special nodes are marked.}
\label{cox}
\vspace{20pt}
\begin{center}
\begin{tabular}{|cc@{\quad}|cc|}
\hline
\raisebox{0pt}{${\bf A_n}$ $(n\ge 1)$}  & 
\raisebox{0pt}{\epsfig{file=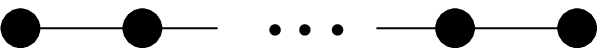,width=0.2\linewidth}}&
\multicolumn{2}{c|}{
\begin{tabular}{cc}
\raisebox{0pt}[20pt][5pt]{${\bf \widetilde A_1}$} & 
\raisebox{0pt}[20pt][5pt]{\epsfig{file=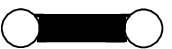,width=0.052\linewidth}}\\
\raisebox{3pt}[15pt][14pt]{${\bf \widetilde A_n}$ $(n\ge 2)$}  & 
\raisebox{-8pt}[25pt][7pt]{\epsfig{file=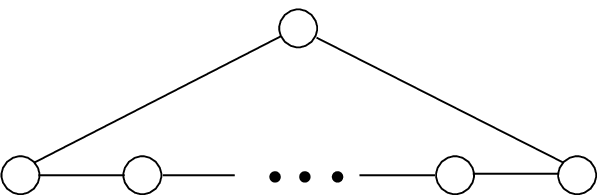,width=0.2\linewidth}}
\end{tabular}
}\\
\hline
\raisebox{-1pt}[23pt][7pt]{${\bf B_n=C_n}$} & 
\raisebox{-7pt}[23pt][7pt]{\epsfig{file=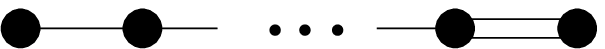,width=0.2\linewidth}}&
\raisebox{7pt}[23pt][7pt]{${\bf \widetilde B_n}$ $(n\ge 3)$}  & 
\raisebox{-0pt}[30pt][7pt]{\epsfig{file=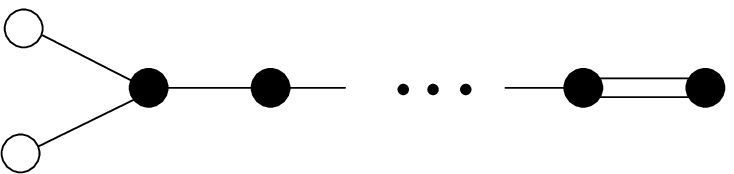,width=0.2\linewidth}}\\
\cline{3-4}
\raisebox{7pt}[15pt][7pt]{ $(n\ge 2)$} & 
&
\raisebox{-0pt}[15pt][7pt]{${\bf \widetilde C_n}$ $(n\ge 2)$}  & 
\raisebox{-0pt}[15pt][7pt]{\epsfig{file=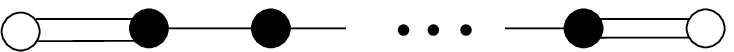,width=0.2\linewidth}}\\
\hline
\raisebox{7pt}[23pt][7pt]{${\bf D_n}$ $(n\ge 4)$} & 
\raisebox{0pt}[30pt][7pt]{\epsfig{file=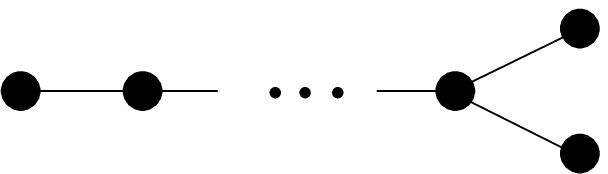,width=0.2\linewidth}}&
\raisebox{7pt}[23pt][7pt]{${\bf \widetilde D_n}$ $(n\ge 4)$} & 
\raisebox{0pt}[30pt][7pt]{\epsfig{file=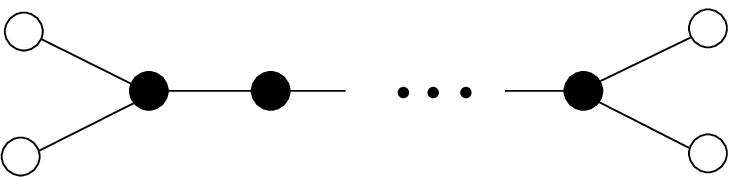,width=0.2\linewidth}}\\
\hline
\raisebox{0pt}[15pt][7pt]{${\bf G_2^{(m)}}$}  & 
\psfrag{m}{{\footnotesize $m$}}
\raisebox{0pt}[15pt][7pt]{\epsfig{file=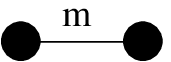,width=0.052\linewidth}}&
\raisebox{0pt}[15pt][7pt]{${\bf \widetilde G_2}$} & 
\raisebox{0pt}[15pt][7pt]{\epsfig{file=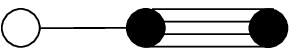,width=0.082\linewidth}}\\
\hline
\raisebox{0pt}[15pt][7pt]{${\bf F_4}$}  & 
\raisebox{0pt}[15pt][7pt]{\epsfig{file=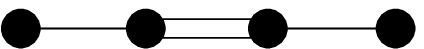,width=0.13\linewidth}}&
\raisebox{-0pt}[15pt][7pt]{${\bf \widetilde F_4}$} & 
\raisebox{-0pt}[15pt][7pt]{\epsfig{file=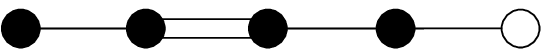,width=0.16\linewidth}}\\
\hline
\raisebox{8pt}[15pt][7pt]{${\bf E_6}$}  & 
\raisebox{0pt}[30pt][7pt]{\epsfig{file=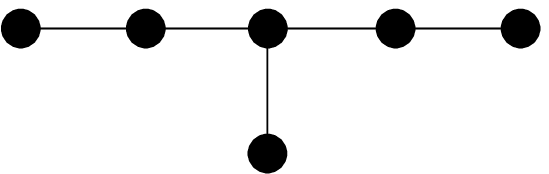,width=0.16\linewidth}}&
\raisebox{8pt}[35pt][7pt]{${\bf \widetilde E_6}$} & 
\raisebox{-8pt}[40pt][17pt]{\epsfig{file=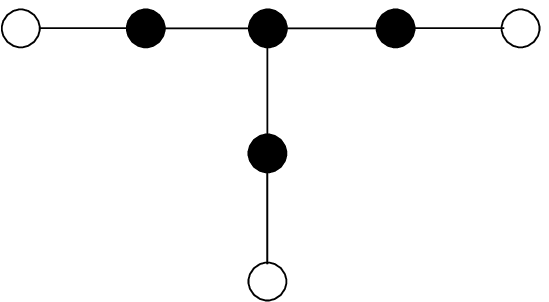,width=0.16\linewidth}}\\
\hline
\raisebox{5pt}[25pt][7pt]{${\bf E_7}$}  & 
\raisebox{-0pt}[30pt][7pt]{\epsfig{file=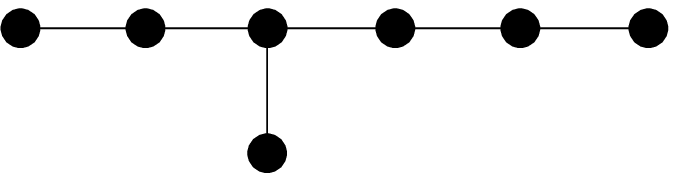,width=0.2\linewidth}}&
\raisebox{5pt}[25pt][7pt]{${\bf \widetilde E_7}$} & 
\raisebox{-0pt}[25pt][7pt]{\epsfig{file=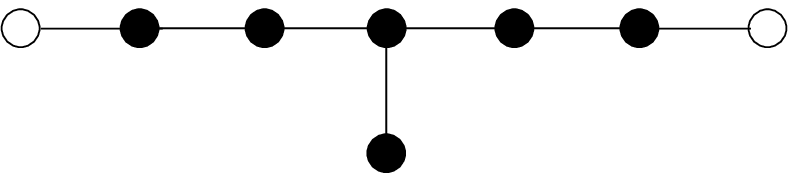,width=0.22\linewidth}}\\
\hline
\raisebox{5pt}[25pt][7pt]{${\bf E_8}$}  & 
\raisebox{-0pt}[30pt][7pt]{\epsfig{file=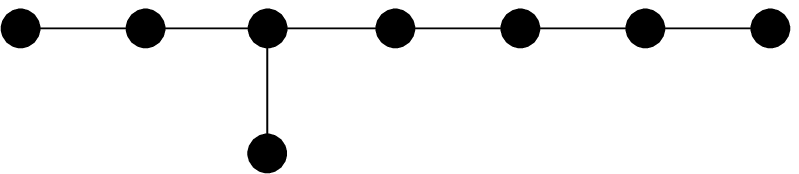,width=0.24\linewidth}}&
\raisebox{5pt}[25pt][7pt]{${\bf \widetilde E_8}$} & 
\raisebox{-0pt}[25pt][7pt]{\epsfig{file=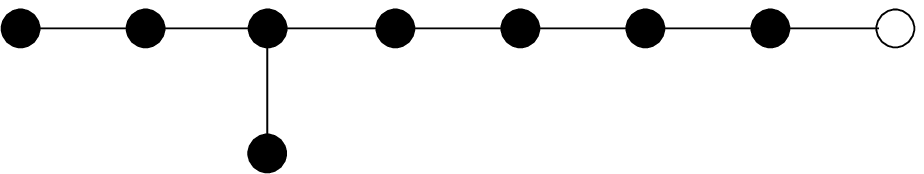,width=0.25\linewidth}}\\
\hline
\raisebox{0pt}[15pt][7pt]{${\bf H_3}$}  & 
\raisebox{0pt}[15pt][7pt]{\epsfig{file=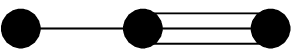,width=0.1\linewidth}}&
& 
\\
\cline{1-2}
\raisebox{0pt}[15pt][7pt]{${\bf H_4}$}  & 
\raisebox{0pt}[15pt][7pt]{\epsfig{file=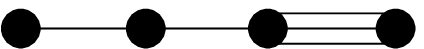,width=0.13\linewidth}}&
& 
\\
\hline

\end{tabular}
\end{center}

\end{table}


%% file: res_e6.txt
515583
128885
104443
104126
104021
104011
64447
40671
39646
 39583
 39573
 35838
 35695
 35629
 35622
 35131
 35128

%% file: res_e7.txt
33554431
33554375
33148285
33148159
16576223
16574174
16574039
16540255
15062365
 15062363
 8288107
 8287019
 7273471
 7273323
 7273255
 7273249
 6813566
 6813547
 6813536
 6813517
 6746927
 6746917
 6746893
 6740393
 6709105
 6709103
 6708972
 6706915
 6706849
 6706529
 6483890
 6483689
 6483675
 6483674
 6482603
 6482587
 6482586
 6480625
 6480411
 6480410
 6466539
 6466515
 6466502
 6462463
 6462076
 6462033
 6462031
 6462029
 6462028
 3112959
 3112885
 3112884
 2579953
 2578873
 2359295
 2359289
 2358204
 2358201
 2358199
 2358196
 2358193
 2353855
 2353852
 2353846
 2353845
 2353843
 2352764
 2352719
 2325497
 2325429
 2325425
 2322165
 2321148
 2321078
 2321074
 2320093
 2320092
 2320044
 2320022
 2320020
 2319949
 2318901
 2318861
 2305919
 2305916
 2305909
 2305885
 2305876
 2305868
 2304831
 2304821
 2304790
 2304780
 2301661
 2301659
 2301656
 2301596
 2301594
 2301593
 2301591
 2301590
 2301588
 2301587
 2301586
 2301584
 2177023
 2177012
 2177006
 2177002
 2176996
 2176993
 2176983
 2176982
 2176958
 2176940
 2176938
 2176934
 2176928
 2176910
 2175854
 2175852
 2175847
 2175844
 2175823
 2175822
 2175820
 2175813
 2175812
 2166783
 2166638
 2166635
 2166634
 2166572
 2166567
 2166563
 2166562
 2166079
 2166074
 2166073
 2166059
 2166057
 2166056

%% file: res_e8.txt
16642998271 \ 
16575889407 \   
16575348735 \   
8036015852 \  
8036015839 \  
8036015609 \  
8027492222 \ 
8027491834 \ 
8023230399 \ 
8023229951 \ 
8023229630 \ 
7834757047 \ 
7834756850 \ 
7834756844 \ 
7834756831 \ 
7834756727 \ 
7834756687 \ 
7834756660 \ 
7834756615 \ 
3758096383 \ 
3758096377 \ 
3749707775 \ 
3749707770 \ 
3749572602 \ 
3749572601 \ 
3749572351 \ 
3749572334 \ 
3749572333 \ 
3749572331 \ 
3749572319 \ 
3745310463 \ 
3745310445 \ 
3745310335 \ 
3745310327 \ 
3745310302 \ 
3743179391 \ 
3743179383 \ 
3743179327 \ 
3743179295 \ 
3611160286 \ 
3611159190 \ 
3611158111 \ 
3611158095 \ 
3611158094 \ 
3611157143 \ 
3611157014 \ 
3611156999 \ 
3548144143 \ 
3548143711 \ 
3548143710 \ 
3543645815 \ 
3543645809 \ 
3543613052 \ 
3543613049 \ 
3543612031 \ 
3543612028 \ 
3543612025 \ 
3543612020 \ 
3543612017 \ 
3543611997 \
3543611995 \ 
3543477119 \ 
3543477079 \ 
3543477023 \ 
3523012220 \  
3523012217 \ 
3523012215 \ 
3523012209 \ 
3523012191 \ 
3523012188 \ 
3523012175 \ 
3523012164 \ 
3522978431 \ 
3522978420 \ 
3522978391 \ 
3522978380 \ 
3522978375 \ 
3522978367 \ 
3522978335 \ 
3522978319 \ 
3522978311 \ 
3522741887 \ 
3522741881 \ 
3522741876 \ 
3522741873 \ 
3522741847 \ 
3522741375 \ 
3522741372 \ 
3522741369 \ 
3514759039 \ 
3514759030 \ 
3514759029 \ 
3514759022 \ 
3514759021 \ 
3514759019 \ 
3514725247 \ 
3514725238 \ 
3514725237 \ 
3514725223 \ 
3514725220 \ 
3514725217 \ 
3514725183 \ 
3514725166 \ 
3514725165 \ 
3514725163 \ 
3514725159 \ 
3514725153 \ 
3514463009 \ 
3514454878 \ 
3514454877 \ 
3514454872 \ 
3514454871 \ 
3514454865 \ 
3514454863 \ 
3514454860 \ 
3514454815 \ 
3514454810 \ 
3514454809 \ 
3514454798 \ 
3514454797 \ 
3514454795 \ 
3514454554 \ 
3514454335 \ 
3514454330 \ 
3514454329 \ 
3514454303 \ 
3514454298 \ 
3514454297 \ 
3514365803 \ 
3514364705 \ 
3514355998 \ 
3514355997 \ 
3514355995 \ 
3514353535 \ 
3514353526 \ 
3514353525 \ 
3514353511 \ 
3514353508 \ 
3514353505 \ 
3514352767 \ 
3514352758 \ 
3514352757 \ 
3514352749 \ 
3514352747 \ 
1606418431 \ 
1606418430 \ 
1606350719 \ 
1606350718 \ 
1606350701 \ 
1606350700 \ 
1604219695 \ 
1405091838 \ 
1405091764 \ 
1405091006 \ 
1405026155 \ 
1405026081 \ 
1405024107 \ 
1404993386 \ 
1404991265 \ 
1404990243 \ 
1404989547 \ 
1404989539 \ 
1404531691 \ 
1404531617 \ 
1404531501 \ 
1404525183 \ 
1404525153 \ 
1404517247 \ 
1404517228 \ 
1404501995 \ 
1404499819 \ 
1404498785 \ 
1404498731 \ 
1404492745 \ 
1404492449 \ 
1404492203 \ 
1404491361 \ 
1404491115 \ 
1404483498 \ 
1402926959 \ 
1402926955 \ 
1402926191 \ 
1402926187 \ 
1402926179 \ 
1402926123 \ 
1402894117 \ 
1402894113 \ 
1402893409 \ 
1402893349 \ 
1402893347 \ 
1402893091 \ 
1402892385 \ 
1402892323 \ 
1402453806 \ 
1402402667 \ 
1402402603 \ 
1402402595 \ 
1402402187 \ 
1402402123 \ 
1402401633 \ 
1402401571 \ 
1402401569 \ 
1402401153 \ 
1402401089 \ 
1402394211 \ 
1402394155 \ 
1402394147 \ 
1402394145 \ 
1402393899 \ 
1402388259 \ 
1402386347 \ 
1402386219 \ 
1402369893 \ 
1402369891 \ 
1402369889 \ 
1402368865 \ 
1402368803 \ 
1402361681 \ 
1402361673 \ 
1402361667 \ 
1402361665 \ 
1402361601 \ 
1402361443 \ 
1402361441 \ 
1402361201 \ 
1402361193 \ 
1402361187 \ 
1402361185 \ 
1402360417 \ 
1402360355 \ 
1402360105 \ 
1402360099 \ 
1402353569 \ 
1402352481 \ 
1375400764 \ 
1375400745 \ 
1375400738 \ 
1375400721 \ 
1375400719 \ 
1375393398 \ 
1375393384 \ 
1375393379 \ 
1375393371 \ 
1375393149 \ 
1375393148 \ 
1375359614 \ 
1375359595 \ 
1375359359 \ 
1375359358 \ 
1375359291 \ 
1375359290 \ 
1371241407 \ 
1371241397 \ 
1371241383 \ 
1371241367 \ 
1371241359 \ 
1371241349 \ 
1371239203 \ 
1371239179 \ 
1371239169 \ 
1371232943 \ 
1371232933 \ 
1371232909 \ 
1371232703 \ 
1371232702 \ 
1371232685 \ 
1371232684 \ 
1371173803 \ 
1371173793 \ 
1371173763 \ 
1371173675 \ 
1371173673 \ 
1371173667 \ 
1371173665 \ 
1371173643 \ 
1371173641 \ 
1371173635 \ 
1371173633 \ 
1371167275 \ 
1371167265 \ 
1371167241 \ 
1371167151 \ 
1371167150 \ 
1371167147 \ 
1371167146 \ 
1371165375 \ 
1371165369 \ 
1371165365 \ 
1371165341 \ 
1371165339 \ 
1371165243 \ 
1371165233 \ 
1371165227 \ 
1371165217 \ 
1371165203 \ 
1371165193 \ 
1371165115 \ 
1371165114 \ 
1371164991 \ 
1371164990 \ 
1371164985 \ 
1371164979 \ 
1371164978 \ 
1371102123 \ 
1371099947 \ 
1371099937 \ 
1371097387 \ 
1371097386 \ 
1371097385 \ 
1370190763 \ 
1370190633 \ 
1370190511 \ 
1370190505 \ 
1370190501 \ 
1370190475 \ 
1370190379 \ 
1370190369 \ 
1370190345 \ 
1370190255 \ 
1370190123 \ 
1370190113 \ 
1370189995 \ 
1370189985 \ 
1370189961 \ 
1369112535 \ 
1369108216 \ 
1369103847 \ 
1369103846 \ 
1369103843 \ 
1369103842 \ 
1369102917 \ 
1369101979 \ 
1369101731 \ 
1369101730 \ 
1369101687 \ 
1369036784 \ 
1369034085 \ 
1369034084 \ 
1369034081 \ 
1369034080 \ 
1368973285 \ 
1368973093 \ 
1368972197 \ 
1368971105 \ 
1368968677 \ 
1368968676 \ 
1368968673 \ 
1368968672 \ 
1368968481 \ 
1368967521 \ 
1368967520 \ 
1368966565 \ 
1368966564 \ 
1368966517 \ 
1368932727 \ 
1368127459 \ 
1368127438 \ 
1368127115 \ 
1368126950 \ 
1368126892 \ 
1368126849 \ 
1368126819 \ 
1368126791 \ 
1368126697 \ 
1368126660 \ 
1368126627 \ 
1368126565 \ 
1368126511 \ 
1358886767 \ 
1358886766 \ 
1358886765 \ 
1358886764 \ 
1358886763 \ 
1358886762 \ 
1358852911 \ 
1358852910 \ 
1358852909 \ 
1358852908 \ 
1358852905 \ 
1358852904 \ 
1358755694 \ 
1358755691 \ 
1358755688 \ 
1358754598 \ 
1358754595 \ 
1358754592 \ 
1358751596 \ 
1358751594 \ 
1358751593 \ 
1358751343 \ 
1358751342 \ 
1358751341 \ 
1358751340 \ 
1358751337 \ 
1358751311 \ 
1358751310 \ 
1358721833 \ 
1358717742 \ 
1358717739 \ 
1358717543 \ 
1358717542 \ 
1358717541 \ 
1358717540 \ 
1358717537 \ 
1358717487 \ 
1358717486 \ 
1358717485 \ 
1358717484 \ 
1358717483 \ 
1358717482 \ 
1358717455 \ 
1358717454 \ 
1357312812 \ 
1357312810 \ 
1357312809 \ 
1357312806 \ 
1357312803 \ 
1357312800 \ 
1357312782 \ 
1357312779 \ 
1357312776 \ 
1357308718 \ 
1357308715 \ 
1357308712 \ 
1357308684 \ 
1357308682 \ 
1357308681 \ 
1357308678 \ 
1357308675 \ 
1357308672 \ 
1357308463 \ 
1357308462 \ 
1357308461 \ 
1357308460 \ 
1357308459 \ 
1357308458 \ 
1357308455 \ 
1357308454 \ 
1357308453 \ 
1357308452 \ 
1357308449 \ 
1357308448 \ 
1357308431 \ 
1357308430 \ 
1357308429 \ 
1357308428 \ 
1357308425 \ 
1357308423 \ 
1357308422 \ 
1357308421 \ 
1357308420 \ 
1357308419 \ 
1357308418 \ 
1357306287 \ 
1357306282 \ 
1354760191 \ 
1354760189 \ 
1354760188 \ 
1354726335 \ 
1354726333 \ 
1354726332 \ 
1354595255 \ 
1354595253 \ 
1354591155 \ 
1354591153 \ 
1354591034 \ 
1354591012 \ 
1354590907 \ 
1354590906 \ 
1354590905 \ 
1354590904 \ 
1354590887 \ 
1354590886 \ 
1354590871 \ 
1354590870 \ 
1354563583 \ 
1354562487 \ 
1354558113 \ 
1354558099 \ 
1354558098 \ 
1354557436 \ 
1354557055 \ 
1354557053 \ 
1354557052 \ 
1354557007 \ 
1354557005 \ 
1354557004 \ 
1354525363 \ 
1354525362 \ 
1354525347 \ 
1354525346 \ 
1354525331 \ 
1354525330 \ 
1354523583 \ 
1354523319 \ 
1354523318 \ 
1354523283 \ 
1354523282 \ 
1354523255 \ 
1354523253 \ 
1354523252 \ 
1354523199 \ 
1354523197 \ 
1354523196 \ 
1354523159 \ 
1354523158 \ 
1354523151 \ 
1354523149 \ 
1354523148 \ 
1353186236 \ 
1353186231 \ 
1353186193 \ 
1353186191 \ 
1353180095 \ 
1353180050 \ 
1353180044 \ 
1353180039 \ 
1353179711 \ 
1353179709 \ 
1353179708 \ 
1353179703 \ 
1353179701 \ 
1353179700 \ 
1353179667 \ 
1353179666 \ 
1353179665 \ 
1353179664 \ 
1353179663 \ 
1353179661 \ 
1353179660 \ 
1353179655 \ 
1353179653 \ 
1353179652 \ 
1353052987 \ 
1353052986 \ 
1353052985 \ 
1353052984 \ 
1353052979 \ 
1353052978 \ 
1353052977 \ 
1353052976 \ 
1353052963 \ 
1353052962 \ 
1353052961 \ 
1353052960 \ 
1353052946 \ 
1353052938 \ 
1353052931 \ 
1353050810 \ 
1353050803 \ 
1353050786 \ 
1353050771 \ 
1353050770 \ 
1353050763 \ 
1353050762 \ 
1353050755 \ 
1353050754 \ 
1353021345 \ 
1353021330 \ 
1353021307 \ 
1353021306 \ 
1353021305 \ 
1353021304 \ 
1353021258 \ 
1353020275 \ 
1353020274 \ 
1353020273 \ 
1353020272 \ 
1353020179 \ 
1353020163 \ 
1353020161 \ 
1353015222 \ 
1353015173 \ 
1353015152 \ 
1353015107 \ 
1353015091 \ 
1353015027 \ 
1353014976 \ 
1353014759 \ 
1353014741 \ 
1353014740 \ 
1353013137 \ 
1353013075 \ 
1353013074 \ 
1353013073 \ 
1353013072 \ 
1353013011 \ 
1353013010 \ 
1353013009 \ 
1353013008 \ 
1353012945 \ 
1353012883 \ 
1353012689 \ 
1353012688 \ 
1353012563 \ 
1353012562 \ 
1353012561 \ 
1353012560 \ 
1353012144 \ 
1353012059 \ 
1353012057 \ 
1353012056 \ 
1353012019 \ 
1353012017 \ 
1353011947 \ 
1353011928 \ 
1353011888 \ 
1353011859 \ 
1353011858 \ 
1353011763 \ 
1353011747 \ 
1353011745 \ 
1353011515 \ 
1353011513 \ 
1353011483 \ 
1353011482 \ 
1353011481 \ 
1353011480 \ 
1353011026 \ 
1353010963 \ 
1353010643 \ 
1353010642 \ 
1353010515 \ 
1353010514 \ 
1353010451 \ 
1353010450 \ 
1353010387 \ 
1352989631 \ 
1352989623 \ 
1352989621 \ 
1352989620 \ 
1352989580 \ 
1352977333 \ 
1352977323 \ 
1352977322 \ 
1352977315 \ 
1352977314 \ 
1352977313 \ 
1352977312 \ 
1352977286 \ 
1352976700 \ 
1352976657 \ 
1352976655 \ 
1352976653 \ 
1352976652 \ 
1352976575 \ 
1352976574 \ 
1352976555 \ 
1352976554 \ 
1352976547 \ 
1352976546 \ 
1352976545 \ 
1352976544 \ 
1352976526 \ 
1352976524
 534773759
 534739895
335544319
335544244
335020031
335011327
333413301
332921844
332914347
332914111
332906383
332889013
332881634
332881633
332881404
332881397
332880545
332880316
332880309
332872644
301693877
301693851
301693848
301693843
301693840
301685435
301685424
301685405
301685181
  299892735
  299892687
  299858932
  299858884
  299630591
  299630555
  299630547
  299630535
  299629500
  299629464
  299629456
  299629444
  299622391
  299622131
  299622119
  299622103
  299621879
  299596788
  299596760
  299596752
  299596748
  299589556
  299589296
  299589284
  299589268
  299588604
  299588568
  299588344
  299588316
  299588084
  299491317
  299491045
  299490228
  299489956
  299489776
  299487221
  299487093
  299486705
  299457532
  299456182
  299455228
  299454388
  299454260
  299453438
  299453418
  299453311
  299453290
  299452732
  299452730
  299452607
  299452606
  298545149
  298545148
  298545077
  298545076
  298545017
  298545016
  298545005
  298545004
  298544945
  298544944
  298544933
  298544932
  298544889
  298544877
  298544861
  298544816
  298544804
  298544788
  298544764
  298544744
  298544716
  298544693
  298544645
  298544568
  298544500
  298544444
  298544424
  298544316
  298544296
  285080502
  285080499
  285080496
  285077229
  285077227
  285077214
  285047733
  285047731
  285047728
  285043639
  285043636
  285043633
  285043428
  285043426
  285043375
  285043366
  285043363
  285043359
  285043350
  284977847
  284977844
  284977831
  284977828
  284977815
  284975861
  284975798
  284975767
  284975646
  284975638
  283999158
  283999157
  283999155
  283999152
  283999015
  283999012
  283999010
  283999009
  283998903
  283998900
  283998898
  283998897
  283998871
  283998868
  283998866
  283998759
  283998756
  283998754
  283998753
  283998726
  283998725
  283998519
  283998517
  283998516
  283998515
  283998513
  283998512
  283997111
  283997108
  283997105
  283996839
  283996836
  283996833
  283996823
  283996820
  283996818
  283996727
  283996724
  283996721
  283996679
  283996676
  283996343
  283996340
  283996338
  283996337
  283634359
  283634356
  283634353
  283634351
  283634346
  283634345
  283634340
  283634338
  283634337
  283634332
  283634330
  283634327
  283634322
  283634311
  283634308
  283634306
  283507637
  283507629
  283507627
  283507620
  283507618
  283507617
  282984447
  282984443
  282983352
  282980349
  282980345
  282980077
  282980073
  282980063
  282951675
  282950588
  282950584
  282949483
  282947508
  282947504
  282947297
  282947293
  282946558
  282946282
  282946276
  282946272
  282946268
  282946223
  282946219
  282946214
  282946210
  282946206
  282946199
  282913524
  282912508
  282912412
  282881789
  282881757
  282880878
  282880702
  282880698
  282880686
  282880682
  282880662
  282879668
  282879580
  282879572
  282878715
  282878685
  282878655
  282878651
  282878614
  282878495
  282878487
  281505533
  281505529
  281505517
  281505513
  281505505
  281505501
  281505497
  281505489
  281505477
  281504500
  281504496
  281504488
  281504484
  281504480
  281504472
  281504468
  281504464
  281504460
  281504446
  281504426
  281504422
  281504418
  281504410
  281504406
  281504402
  281504390
  281504386
  281503484
  281503480
  281503472
  281503468
  281503464
  281503444
  281503436
  281503432
  281503416
  281503412
  281503408
  281503396
  281503392
  281503388
  281503368
  281503364
  281503360
  281411579
  281411562
  281411555
  281410544
  281410536
  281410529
  281410492
  281410488
  281410477
  281410473
  281410468
  281410464
  281408509
  281408500
  281408485
  281408445
  281408441
  281408436
  281408432
  281408428
  281408424
  281408285
  281408281
  281408276
  281408272
  281408261
  281408257
  281407453
  281407436
  281407429
  281407380
  281407372
  281407365
  281407261
  281407252
  281407237
  281407197
  281407188
  281407173
  281405401
  281405336
  281405329
  281405245
  281405236
  281405228
  281405140
  281405132
  281405125
  281405020
  281405013
  281404996
  281404221
  281404212
  281404197
  281403967
  281403950
  281403943
  281403935
  281403926
  281403911
  281402329
  281402320
  281402312
  281402269
  281402265
  281402260
  281402256
  281402252
  281402248
  281402141
  281402137
  281402132
  281402128
  281402124
  281402120
  281402077
  281402068
  281402060
  281401945
  281401936
  281401928
  281401885
  281401876
  281401868
  281401757
  281401753
  281401149
  281401140
  281401044
  281400957
  281400948
  281400940
  281400895
  281400886
  281400878
  281400690
  281400671
  281399293
  281399228
  281399221
  281399101
  281399092
  281399077
  281399005
  281398876
  281398861
  281398813
  281398796
  281398789
  281398513
  281398505
  281378815
  281378799
  281378795
  281378785
  281377780
  281377770
  281377764
  281377760
  281377702
  281377698
  281376766
  281376756
  281376750
  281376746
  281376702
  281376692
  281376676
  281376672
  281376598
  281376582
  281369564
  281369560
  281369554
  281369548
  281369544
  281369496
  281369494
  281369490
  281369480
  281369436
  281369432
  281369430
  281369420
  281369416
  281369238
  281369228
  281368988
  281368984
  281368982
  281368534
  281368530
  281368524
  281368518
  281368514
  281368468
  281368464
  281368462
  281368452
  281368448
  281368446
  281368436
  281368432
  281368430
  281368214
  281368210
  281368204
  281368198
  281368194
  281368182
  281368172
  281367960
  281367958
  281367954
  281367926
  281367922
  281367912
  281367902
  281367892
  281367888
  281367422
  281367418
  281367412
  281367196
  281367166
  281367162
  281367156
  281367140
  281367103
  281367093
  281367077
  281367073
  281366900
  281366896
  281366876
  281366833
  281366813
  281366809
  281366807
  281366527
  281366517
  281366501
  281366497
  281366462
  281366452
  281366436
  281366432
  281366398
  281366388
  281366372
  281366368
  281365912
  281365910
  281365906
  281365848
  281365846
  281365842
  281365785
  281365783
  281365779
  281365749
  281365739
  281365733
  281365729
  281365502
  281365430
  281365374
  281365364
  281365358
  281365354
  281364724
  281364714
  281364708
  281364704
  281364652
  281364642
  276824063
  276824062
  276790199
  276790198
  276758527
  276757429
  276756350
  276756335
  276756334
  276756331
  276724662
  276722548
  276722487
  276722477
  276722476
  276722471
  276722470
  276722467
  275676085
  275676084
  275675957
  275675956
  275675941
  275675940
  275675937
  275675936
  275217333
  275217332
  275217313
  275217312
  275217285
  275217284
  275215285
  275215284
  275215265
  275215264
  275215237
  275215236
  275215157
  275215156
  275215137
  275215136
  275215109
  275215108
  275209715
  275209141
  275209140
  275209009
  275209008
  275208997
  275208996
  275208925
  275208869
  275208868
  275208853
  275208852
  275208632
  275176100
  275176076
  275175870
  275175860
  275175855
  275175854
  275175850
  275151781
  275151780
  275151749
  275151748
  275149733
  275149732
  275149701
  275149700
  275149605
  275149604
  275149573
  275149572
  275144597
  275144452
  275144183
  275144167
  275144166
  275144162
  275143605
  275143604
  275143585
  275143584
  275143573
  275143572
  275143473
  275143472
  275143461
  275143460
  275143441
  275143440
  275143429
  275143428
  275143333
  275143332
  275143301
  275143300
  275143156
  275143141
  275143140
  275143137
  275143101
  275143085
  275143084
  275143080
  275142130
  275142119
  275142115
  275142114
  275142064
  275142053
  275142049
  275142048
  275142006
  275142003
  275142002
  275141987
  275141557
  275141556
  275141425
  275141424
  275141413
  275141412
  275141301
  275141300
  275141281
  275141280
  275141269
  275141268
  275141208
  275141157
  275141156
  275141125
  275141124
  275141049
  275140980
  275140977
  275140976
  275140924
  275140921
  275140920
  275140905
  274987958
  274987957
  274987955
  274987952
  274985783
  274985780
  274985778
  274985777
  274985767
  274985766
  274985765
  274985764
  274985763
  274985762
  274985761
  274985760
  274983713
  274983607
  274983604
  274983602
  274983601
  274983575
  274983574
  274981681
  274981680
  274981543
  274981540
  274981538
  274981537
  274981430
  274981429
  274981427
  274981424
  274981415
  274981414
  274981413
  274981412
  274981411
  274981410
  274981409
  274981408
  274981399
  274981398
  274922423
  274922421
  274922420
  274918324
  274918320
  274918067
  274918066
  274918065
  274918064
  274918039
  274918038
  274918035
  274918034
  274916003
  274916001
  274916000
  274915895
  274915893
  274915892
  274915863
  274915862
  274915859
  274915858
  274661375
  274659198
  274659181
  274659179
  274626414
  274625388
  274625327
  274625315
  273580031
  273579884
  273579882
  273579839
  273579815
  273579809
  273154047
  273154026
  273154017
  273153991
  273150975
  273150954
  273150945
  273150919
  273150783
  273150762
  273150753
  273150727
  273145855
  273145657
  273145639
  273145575
  273145567
  273087406
  273087396
  273087366
  273086382
  273086372
  273086342
  273086318
  273086308
  273086278
  273080252
  273080226
  273080222
  273080186
  273080164
  273080152
  273080146
  273080140
  273080039
  273080015
  273079789
  273079208
  273079152
  273079150
  273079024
  273079022
  273078956
  273078926
  273078718
  273078698
  273078198
  273078128
  273078126
  273078006
  273077992
  273077982
  273077940
  273077930
  273077916
  273077910
  273077868
  273077838
  273077624
  273077618
  273055724
  273055719
  273055687
  273052652
  273052647
  273052615
  273052460
  273052455
  273052423
  273047551
  273047521
  273047519
  273047498
  273047353
  273047335
  273047321
  273047308
  273047303
  273047271
  273047244
  273047239
  273047023
  273047010
  273045492
  273045482
  273045471
  273045460
  273045450
  273045298
  273045292
  273045273
  273045266
  273045260
  273045255
  273045238
  273045224
  273045213
  273045206
  273045192
  273045162
  273045151
  273045130
  273045100
  273045068
  273045063
  273044978
  273044970
  273044967
  273044961
  273044904
  273044899
  273044859
  273044854
  273044848
  273044835
  273044479
  273044281
  273044263
  273044223
  273044193
  273044191
  273044090
  273044068
  273044058
  273044049
  273044047
  273044007
  273043980
  273043975
  273043833
  273043826
  273043770
  273043753
  272826367
  272826365
  272825278
  272825276
  272824186
  272824184
  272823222
  272823220
  272823218
  272823216
  272823154
  272823152
  272823099
  272823097
  272823079
  272823077
  272823075
  272823073
  272822271
  272822265
  272822198
  272822192
  272822126
  272822120
  272822055
  272822049
  272822011
  272822009
  272821981
  272821979
  272821180
  272821178
  272821036
  272821034
  272820978
  272820976
  272820922
  272820920
  272820902
  272820900
  272820898
  272820896
  272820890
  272820092
  272820090
  272820020
  272820018
  272819946
  272819944
  272819874
  272819872
  272819838
  272819836
  272819804
  272819802
  272819062
  272819056
  272819007
  272819001
  272818915
  272818913
  272818870
  272818868
  272818866
  272818864
  272818858
  272818856
  272818806
  272818804
  272818772
  272818770
  272818751
  272818749
  272818727
  272818725
  272818723
  272818721
  272818717
  272818715
  272807868
  272807866
  272807862
  272807856
  272807854
  272807848
  272807844
  272807842
  272807838
  272807832
  272807828
  272807826
  272807820
  272807818
  272807814
  272807808
  272804730
  272804728
  272804722
  272804720
  272804714
  272804712
  272804706
  272804704
  272804698
  272804696
  272804690
  272804688
  272804682
  272804680
  272804674
  272804672
  272803774
  272803766
  272803764
  272803754
  272803746
  272803744
  272803726
  272803724
  272803718
  272803716
  272800622
  272800620
  272800618
  272800616
  272800372
  272800368
  272800366
  272800364
  272800362
  272800360
  272800334
  272800332
  272800330
  272800328
  272789503
  272789497
  272789229
  272789224
  272789211
  272787196
  272787193
  272787114
  272787113
  272787063
  272787060
  272787036
  272787034
  272786155
  272786152
  272786047
  272786044
  272786002
  272785947
  272771071
  272771068
  272771063
  272771060
  272771034
  272771033
  272771026
  272771025
  272771023
  272771020
  272771015
  272771012
  272770746
  272770732
  272770730
  272770722
  272770714
  272770700
  272770698
  272770690
  272769023
  272769017
  272768639
  272768636
  272768633
  272768599
  272768594
  272768593
  272768591
  272768588
  272768585
  272768575
  272768569
  272768527
  272768524
  272768521

%% file: res_f4.txt
32560
4100
3232
3220
1213
1159
1153
1124
1064
 827
 825